\documentclass{amsart}

 \theoremstyle{definition}
 
 \theoremstyle{remark}
 
 \numberwithin{equation}{subsection}
\begin{document}

\title[Jordan counterparts of Rickart and Baer $*$-algebras, II]
{Jordan counterparts of Rickart and Baer $*$-algebras, II}

\author{Shavkat Ayupov}

\address{Institute of Mathematics, National University of Uzbekistan, Tashkent, Uzbekistan}

\email{sh$_-$ayupov@mail.ru}

\medskip

\author{Farhodjon Arzikulov}

\address{Department of Mathematics and Physics,
Andizhan State University, Andizhan, Uzbekistan.}

\email{arzikulovfn@rambler.ru}

\thanks{This paper partially supported by TWAS, The Abdus Salam, International
Centre, for Theoretical Physics (ICTP), Grant:
13-244RG/MATHS/AS$_-$I-UNESCO FR:3240277696}

\begin{abstract}
We introduce and investigate  new classes of Jordan algebras
which are close to but wider than Rickart and Baer Jordan algebras considered in our previous paper.
Such Jordan algebras are called RJ- and BJ-algebras respectively.
Criterions are given for a Jordan algebra to be a BJ-algebra. Also, it is proved that every finite dimensional
Jordan algebra without nilpotent elements, which have square roots, is a BJ-algebra.
\end{abstract}

\maketitle


\section*{Introduction}

In the previous paper \cite{AA} the authors  introduced and investigated Rickart and Baer Jordan algebras.
In \cite{Arz} it is proved that most of known examples of Jordan Banach algebras
satisfy the following algebraic condition: for every subset $S$ of
a Jordan Banach algebra $A$, the Jordan annihilator $P^\perp :=\{a\in
A:U_ax=0, (\forall x\in S)\}$ of the set $P=\{a^2: a\in S\}$ is an
inner ideal, generated by an idempotent $e\in A$, i.e.
$P^\perp=U_e(A)$. In \cite{AA} a Jordan algebra satisfying this condition was called
a Baer Jordan algebra.

The present paper is devoted to investigating of Jordan algebras with
more general conditions than the conditions for Rickart and Baer Jordan algebras introduced in \cite{AA}.
Such algebras are called RJ- and
BJ-algebras respectively. By the Jordan-von Neumann-Wigner theorem \cite{JNW}
every finite-dimensional, formally real, unital Jordan algebra
over the field of real numbers is a BJ- and RJ-algebra
\cite{HOS}. The definitions of RJ- and BJ-algebras
are quite compatible with the definitions of Rickart and Baer $*$-algebras. It
turns out that RJ- and BJ-algebras have no nilpotent elements, which have  square roots. Also,
a Jordan algebra $A$ is a BJ-algebra if and only if the quotient $A/\mathcal{D}eg(A)$
is a BJ-algebra, where  $\mathcal{D}eg(A)$ is the degenerate radical.

Finally, we prove that every finite dimensional
Jordan algebra without nilpotent elements, which have square roots, is a BJ-algebra.

We use \cite{AA}, \cite{KMcC} as a standard reference on
notations and terminology. Throughout the paper we consider Jordan
algebras $A$ over a field $F$ of characteristic $\neq 2$.

\bigskip

\section{RJ-algebras}

Let $A$ be a Jordan algebra. Introduce the following notations:
$$
A^2=\{a^2: a\in A\}, \{abc\}=(ab)c+(cb)a-(ac)b,
$$
$$
U_ab=2(ab)a-a^2b, U_a(A)=\{U_ax: x\in A\},
$$
$$
S^\perp :=\{a\in A:U_ax=0, (\forall x\in S)\},
^\perp S:=\{x\in A:U_ax=0, (\forall a\in S)\}
$$
and
$$
^\perp S\cap A^2=^\perp S\cap A^2, {S}^\perp\cap A^2={S}^\perp\cap A^2, U_a(A)\cap A^2=U_a(A)\cap A^2.
$$

For the Jordan algebra $A$ consider the following condition:

(A) for every element $x\in A$ there exists an idempotent $e\in
A$ such that

 $$
 ^\perp \{x\}\cap A^2=U_e(A)\cap A^2.
$$
{\it Definition.} A Jordan algebra satisfying the condition (A) is
called {\it an RJ-algebra}.

It should be noted that initially analogues of RJ-algebras have been mentioned in
\cite{KMP}.

If for every element $x\in A^2$ of a Jordan algebra $A$ there exists an idempotent $e\in
A$ such that
$$
\{x\}^\perp =U_e(A),
$$

then $A$ is called a Rickart Jordan algebra \cite{AA}.

{\bf Proposition 1.1.}
{\it Every Rickart Jordan algebra is an RJ-algebra.}

{\it Proof.} This proposition follows by theorem 1.6 in \cite{AA}. $\triangleright$

For a (real) $*$-algebra $M$ , and a nonempty subset $S$ in $M$, we
set
$$
R(S)=\{x\in M: sx=0 (\forall s\in S)\},
$$
and call $L(S)$ the right-annihilator of $S$.

Let $M$ be a (real) $*$-algebra and $M_{sa}=\{a\in M: a^*=a\}$.
Then the set $M_{sa}$ is a Jordan
algebra with respect to the multiplication $a\cdot b=1/2(ab+ba)$.

{\it Definition.} A (real) Rickart  *-algebra is a (real) *-algebra $M$ such
that, for each $x\in M$, $R(\{x\})=gM$ for an appropriate projection $g$ in $M$
(i. e. $g^*=g$, $g^2=g$).

{\bf Proposition 1.2.}
{\it Let $M$ be a (real) Rickart  $*$-algebra. Then
$M_{sa}$ is an RJ-algebra.}

{\bf Proof.} The Jordan algebra $A = M_{sa}$ is an RJ-algebra
by proposition 1.1 in \cite{AA}. $\triangleright$

The following Jordan algebras are additional examples of RJ-algebras:

1.  Let
$H_n(F)$ be the Jordan algebra of $n\times n$ dimensional self-adjoint matrices over $F$,
where $F={\Bbb R}$ (reals), ${\Bbb C}$ (complexes), ${\Bbb H}$ (quaternions) or
${\Bbb Q}$ (octonions). By proposition 1.2 and proposition 3.4 in \cite{AA} $H_n(F)$ is an RJ-algebra.
So the set $\mathcal{R}$ of all infinite sequences
with components from $H_n(F)$ and with finite number of
nonzero elements is also an RJ-algebra with respect to the
componentwise algebraic operations.

Recall that the {\it unital hull} of a Jordan algebra $A$ is the following set
$$
\hat{A}:=F\hat{1}\oplus A
$$
with the Jordan multiplication $(\alpha\hat{1}\oplus x)(\beta\hat{1}\oplus y):=\alpha\beta\hat{1}\oplus(\alpha y+\beta x + xy)$.
Every Jordan algebra $A$ can be imbedded as an ideal into its unital hull. An element $z\in A$  is called
{\it trivial} (see \cite{KMcC}) if its U-operator is trivial on the unital hull $\hat{A}$, i.e.
$U_z(\hat{A})=0$, that is equivalent to $U_z(A)=0$, $z^2=0$ \cite{KMcC}.

2. Set
$$
{\bf 0}=\left(
           \begin{array}{cc}
             0 & 0 \\
             0 & 0 \\
           \end{array}
         \right),
{\bf 1}=\left(
           \begin{array}{cc}
             1 & 0 \\
             0 & 1 \\
           \end{array}
         \right),
e_{1,2}=\left(
         \begin{array}{cc}
           0 & 1 \\
           0 & 0 \\
         \end{array}
       \right).
$$
Then $\mathcal{A}={\Bbb R}{\bf 1}+{\Bbb R}e_{1,2}$ is an associative algebra.
$\mathcal{A}$ is a Jordan algebra with respect to the multiplication $x\cdot y=\frac{1}{2}(xy+yx)$.
We have
$$
^\perp\{{\bf 0}\}\cap \mathcal{A}^2=U_{\bf 1}(\mathcal{A})\cap \mathcal{A}^2,
^\perp\{{\bf 1}\}\cap \mathcal{A}^2=U_{\bf 0}(\mathcal{A})\cap \mathcal{A}^2,
^\perp\{e_{1,2}\}\cap \mathcal{A}^2=U_{\bf 1}(\mathcal{A})\cap \mathcal{A}^2.
$$
Also for every nonzero number $\lambda$ in ${\Bbb R}$ we have
$$
^\perp\{\lambda{\bf 1}\}\cap \mathcal{A}^2=U_{\bf 0}(\mathcal{A})\cap \mathcal{A}^2,
^\perp\{\lambda e_{1,2}\}\cap \mathcal{A}^2=U_{\bf 1}(\mathcal{A})\cap \mathcal{A}^2.
$$
Now, we take nonzero numbers $\lambda$ and $\mu$ in ${\Bbb R}$. For arbitrary element
$\alpha{\bf 1}+\beta e_{1,2}$ in $\mathcal{A}$ we have
$$
U_{\lambda{\bf 1}+\mu e_{1,2}}(\alpha{\bf 1}+\beta e_{1,2})=\lambda^2\alpha{\bf 1}+\lambda^2\beta e_{1,2}
+2\lambda\mu\alpha e_{1,2}\neq 0
$$
if $(\alpha,\beta)\neq (0,0)$.
Hence
$$
^\perp\{\lambda{\bf 1}+\mu e_{1,2}\}\cap \mathcal{A}=U_{\bf 0}(\mathcal{A})\cap \mathcal{A}^2.
$$
Thus $\mathcal{A}$ is an RJ-algebra. Note that $\mathcal{A}$ is an example of
an RJ-algebra which has a trivial element $e_{1,2}$  and nonzero nilpotent radical $\mathcal{N}il(\mathcal{A})={\Bbb R}e_{1,2}$.

3. Let $M_n({\Bbb R})$ be the
matrix algebra over $\Bbb R$, of $n>1$ dimensional matrices of the form
$$
\left[%
\begin{array}{cccc}
a^{1,1} & a^{1,2} & \cdots & a^{1,n}\\
a^{2,1} & a^{2,2} & \cdots & a^{2,n}\\
\vdots & \vdots & \ddots & \vdots\\
a^{n,1} & a^{n,2} & \cdots & a^{n,n}\\
\end{array}%
\right],
a^{i,j}\in \Re, i,j=1,2,\dots ,n.
$$
Let $\{e_{i,j}\}_{i,j=1}^n$ be the set of matrix units in
$M_n({\Bbb R})$, i.e. $e_{i,j}$ is a matrix with components
$a^{i,j}=1$ and $a^{k,l}=0$ if $(i,j)\neq(k,l)$ and
a matrix $a\in M_n({\Bbb R})$ is written as $a=\sum_{k,l=1}^n
a^{k,l}e_{k,l}$, where $a^{k,l}\in {\Bbb R}$ for $k,l=1,2,\dots, n$.

Let ${\bf 0}$ be the zero matrix and ${\bf 1}$ be the identity matrix in
$M_n({\Bbb R})$. Suppose that  $n$ is even and $n=2k$ for some natural number $k$.

Then
$$
\mathcal{A}={\Bbb R}{\bf 1}+\sum_{i=1}^k {\Bbb R}e_{2i-1,2i}
$$
is an associative algebra and it becomes a Jordan algebra when equipped with the
multiplication $x\cdot y=\frac{1}{2}(xy+yx)$.
We have
$$
^\perp\{{\bf 0}\}\cap \mathcal{A}^2=U_{\bf 1}(\mathcal{A})\cap \mathcal{A}^2,
^\perp\{{\bf 1}\}\cap \mathcal{A}^2=U_{\bf 0}(\mathcal{A})\cap \mathcal{A}^2,
$$
$$
^\perp\{e_{2i-1,2i}\}\cap \mathcal{A}^2=U_{\bf 1}(\mathcal{A})\cap \mathcal{A}^2, i=1,2,\dots,k.
$$
Also for arbitrary numbers $\lambda\neq 0$, $\lambda_1$, $\lambda_2$, $\dots$, $\lambda_k$ in ${\Bbb R}$ we have
$$
^\perp\{\lambda{\bf 1}\}\cap \mathcal{A}^2=U_{\bf 0}(\mathcal{A})\cap \mathcal{A}^2,
^\perp\{\sum_{i=1}^k \lambda_i e_{2i-1,2i}\}\cap \mathcal{A}^2=U_{\bf 1}(\mathcal{A})\cap \mathcal{A}^2.
$$
Now, we take an arbitrary element
$\lambda {\bf 1}+\sum_{i=1}^k \lambda_i e_{2i-1,2i}$ in $\mathcal{A}$ with a nonzero number $\lambda$ in ${\Bbb R}$.
Then
$$
U_{\lambda {\bf 1}+\sum_{i=1}^k \lambda_i e_{2i-1,2i}}(\alpha{\bf 1}+\sum_{i=1}^k \alpha_i e_{2i-1,2i})=
$$
$$
\lambda^2\alpha{\bf 1}+2\lambda\alpha(\sum_{i=1}^k \lambda_i e_{2i-1,2i})+\lambda^2\sum_{i=1}^k \alpha_i e_{2i-1,2i}\neq 0
$$
if $(\alpha,\alpha_1,\dots, \alpha_k)\neq (0,0,\dots,0)$.
Hence
$$
^\perp\{\lambda {\bf 1}+\sum_{i=1}^k \lambda_i e_{2i-1,2i}\}\cap \mathcal{A}=U_{\bf 0}(\mathcal{A})\cap \mathcal{A}^2.
$$
So $\mathcal{A}$ is an RJ-algebra. Note that this is an example of
an RJ-algebra which has a trivial element and nonzero nilpotent radical $\mathcal{N}il(\mathcal{A})=\sum_{i=1}^k {\Bbb R}e_{2i-1,2i}$.

{\bf Lemma 1.3.}
1) {\it Let $A$ be an RJ-algebra. Then

a) $A$ has an element $e$ such that $ae=ea=a$ for every $a\in A^2$,

b) the algebra $A$ has no nilpotent element, with a square root.}

{\bf Proof.} a). Take $x=0$. We have $^\perp\{x\}\cap A^2=U_e(A)\cap A^2$ for
some idempotent $e$ in $A$. But $^\perp\{x\}\cap A^2=A\cap A^2=A^2$. Hence $U_e(A)\cap A^2=A^2$.
Then $ea=eU_eb=U_eb=a$ and $ae=(U_eb)e=U_eb=a$ for all $a\in A^2$.

b). Fix the above idempotent $e$. Take an element $a\in A^2$. Suppose $a^2=0$,
 then $U_ae=0$ and $a\in {^\perp\{e\}}$. By the condition
there exists an idempotent $f\in A$ such that $^\perp\{e\}\cap A^2
=U_f(A)\cap A^2$. Since $U_ff=f$ we have $f\in U_f(A)$ and $f\in
^\perp\{e\}\cap A^2$. Hence $U_ef=0$, i.e. $f=0$. Therefore $^\perp\{e\}\cap A^2
=\{0\}$ and $a=0$. Suppose $a^3=0$, $a\in A$; then $a^4=0$ and by
the previous part of the proof $a^2=0$ since $a^4=(a^2)^2=0$.
Hence $a=0$.

 Further we shall apply the induction.Fix a natural number $n$.
 Consider  $a\in A^2$. Suppose that  $a^k=0$,  $k\leq n,$ implies that $a=0.$
Let $a^{n+1}=0$. Then, if $n=2m$ then $a^{2m+2}=(a^{m+1})^2=0$. By
the previous part of the proof $a^{m+1}=0$ and by inductive
assumption $a=0$. if $n=2m+1$ then
$a^{n+1}=a^{2m+2}=(a^{m+1})^2=0$. Again by the previous part of
the proof $a^{m+1}=0$ and $a=0$. Hence by induction we obtain that
for  $a\in A^2$, if $a^n=0$ for some natural number $n$, then $a=0$ .
The proof is complete.
$\triangleright$

{\it Remark.} 1. Here we give an example of an RJ-algebra without a unit element.
We take the associative algebra $\mathcal{A}=\sum_{i=1}^k {\Bbb R}e_{2i-1,2i}$.
$\mathcal{A}$ is a Jordan algebra with the multiplication $x\cdot y=\frac{1}{2}(xy+yx)$.
But this Jordan algebra has no nonzero idempotent elements and
$$
^\perp\{{\bf 0}\}=U_{\bf 0}(A)\cap A^2.
$$

2. The following Jordan algebras are not RJ-algebras: let
$M_n(F)$ be the $*$-algebra of $n\times n$ dimensional matrices over $F$,
where $F={\Bbb R}$, ${\Bbb C}$, ${\Bbb H}$.
They are Jordan algebras with respect to the multiplication $x\cdot y=\frac{1}{2}(xy+yx)$.
For $n>2$ these Jordan algebras are not RJ-algebras. Indeed, by
lemma 1.3 an RJ-algebra has no nilpotent elements, which have square roots. At
the same time
$$
\left(
  \begin{array}{ccc}
    0 & 1 & 1 \\
    0 & 0 & 1 \\
    0 & 0 & 0 \\
  \end{array}
\right)^4
=
\left(
  \begin{array}{ccc}
    0 & 0 & 1 \\
    0 & 0 & 0 \\
    0 & 0 & 0 \\
  \end{array}
\right)^2=
\left(
  \begin{array}{ccc}
    0 & 0 & 0 \\
    0 & 0 & 0 \\
    0 & 0 & 0 \\
  \end{array}
\right).
$$
Hence the Jordan algebra $M_3(F)$ is not an RJ-algebra by b) of lemma 1.3.
Similarly for each $n>2$ $M_n(F)$ is not an RJ-algebra.

The element $e$ from b) of lemma 1.3  we denote by $1$.

{\bf Lemma 1.4.}
{\it Let $A$ be an RJ-algebra, and let $e$ be an idempotent in $A$. Then
${^\perp\{e\}}\cap A^2=U_{1-e}(A)\cap A^2$.}

{\it Proof.} We have $U_e(1-e)=0$ and $1-e\in {^\perp\{e\}}\cap A^2$. Since $A$ is an RJ-algebra
we have ${^\perp\{e\}}\cap A^2=U_f(A)\cap A^2$ for an idempotent $f\in A$. At the same time,
$U_eU_{1-e}(A)=U_{e(1-e)}(A)=\{0\}$. Hence $U_{1-e}(A)\subseteq U_f(A)$ and $f(1-e)=1-e$.

Since $U_ef=0$ and $f=U_ef+\{(1-e)fe\}+U_{1-e}f$ we have
$f=\{(1-e)fe\}+U_{1-e}f$. Hence
$f=\{(1-e)fe\}+U_{1-e}f=\{(1-e)fe\}+(1-e)$ by $f(1-e)=1-e$.
Also we have $fe=e\{(1-e)fe\}+e(1-e)=\frac{1}{2}\{(1-e)fe\}$ and
$1-e=f(1-e)=\{(1-e)fe\}(1-e)+(1-e)^2=\frac{1}{2}\{(1-e)fe\}+1-e$, $\{(1-e)fe\}=0$.
Therefore $fe=0$ and $f(1-e)=f$.

Thus $f=1-e$. The proof is complete. $\triangleright$

{\it Definition.} Let $A$ be an RJ-algebra, $x\in A$,
and let $^\perp\{x\}\cap A^2=U_e(A)\cap A^2$ for an appropriate idempotent $e$. This
idempotent is unique. Indeed, if $U_e(A)=U_f(A)$ with $e$ and $f$
idempotents, then $fe=e$ and $ef=f$, so $e=fe=ef=f$. Thus the
idempotent in the definition of RJ-algebras is
unique.

For idempotents $e$, $f$ in a Jordan algebra $A$, one writes
$e\leq f$ if $e\in U_f(A)$ that is, $ef=fe=e$. For idempotents
$e$, $f$ in a Jordan algebra $A$, the following conditions are
equivalent: $e\leq f$, $e=ef$, $U_e(A)\subseteq U_f(A)$.

{\it Definition.} A Jordan algebra $A$ is said to be
{\it non-degenerate} if it has no nonzero trivial elements. A Jordan algebra $A$ is said to be
{\it quadratic non-degenerate} if it has no nonzero trivial elements with square roots.

{\bf Theorem 1.5.}
{\it Every RJ-algebra is quadratic non-degenerate.}

{\it Proof.} Let $A$ be an RJ-algebra and let $a$ be a nonzero element in $A^2$. Then
there exists an idempotent $e$ such that
$$
^\perp \{a\}\cap A^2=U_e(A)\cap A^2.
$$
It is clear that $e\neq 1$, where $1$ is the unit element of $A$. Otherwise
$^\perp \{a\}\cap A^2=A^2$ and $U_a1=a^2=0$. Hence $a=0$ that
is impossible by lemma 1.3. So $e<1$. Therefore $U_a(1-e)\neq 0$, i.e. $U_a(A)\neq 0$.
The proof is complete. $\triangleright$

Degeneracy in a Jordan algebra $A$ can be localized in the degenerate
radical $\mathcal{D}eg(A)$, the smallest ideal, the quotient respect to which is nondegenerate \cite{KMP}.

{\bf Theorem 1.6.}
{\it Let $A$ be a Jordan algebra. Then, if the quotient
$A/\mathcal{D}eg(A)$ is an RJ-algebra, then $A$ is an RJ-algebra.}

{\it Proof.}
Let $B=A/\mathcal{D}eg(A)$ and $a$ be an arbitrary nontrivial element in $A$. Then there exists
an idempotent $e$ in $A$ such that
$$
^\perp\{a\}\cap B^2=U_e(B)\cap B^2.
$$
Hence $^\perp\{a\}\cap A^2=U_e(A)\subseteq A^2$. Take an element
$b\in \mathcal{D}eg(A)$. We have, if $b=c^2$ for some element $c\in A$,
then $b^2=(c^2)^2=c^4=0$ and, hence $c=0$. Therefore $b=0$. So
$$
\mathcal{D}eg(A)\cap A^2=\{0\} \,\,\,\text{and}\,\,\, B^2\equiv A^2.
$$
Thus $^\perp\{a\}\cap A^2\equiv (^\perp\{a\})\cap B^2=U_e(B)\cap B^2\equiv U_e(A)\cap B^2=U_e(A)\cap A^2$.
Hence $^\perp\{a\}\cap A^2=U_e(A)\cap A^2$.
Also, for every trivial element $a\in A$ we have $^\perp\{a\}\cap A^2=A^2=U_1(A)\cap A^2$.
So $A$ is an RJ-algebra. The proof is complete. $\triangleright$

{\it Definition.} (1) An element $x$ of an arbitrary unital Jordan algebra $A$
is quasi-invertible if $1-z$ is invertible in $A$; if
$(1-z)^{-1}=1-w$ then $w$ is called the quasi-inverse of $z$ and denoted by $qi(z)$.
The set of quasi-invertible elements of $J$ is denoted by $\mathcal{QI}(A)$. An algebra or
ideal is called quasi-invertible if all its elements are quasi-invertible.

(2) The Jacobson radical $\mathcal{R}ad(A)$ is the maximal quasi-invertible ideal, i.e., the maximal
ideal contained inside the set $\mathcal{QI}(A)$.

{\bf Corollary 1.7.}
{\it Every finite dimensional Jordan algebra without nilpotent elements
having square roots, is an RJ-algebra.}

{\it Proof.} Let $A$ be a finite dimensional unital Jordan algebra.
By radical equality theorem 1.7.2 in \cite{KMcC} $\mathcal{D}eg(A)$
coincides with $\mathcal{R}ad(A)$, i.e. $\mathcal{D}eg(A)=\mathcal{R}ad(A)$.
Hence by enlightenment structure theorem in page 79 of \cite{KMcC}, by theorems 1.5, 3.4 in \cite{AA} and
by proposition 1.2 the quotient $A/\mathcal{D}eg(A)$ is an RJ-algebra. Therefore
by theorem 1.6 $A$ is an RJ-algebra. The proof is complete. $\triangleright$.

{\bf Proposition 1.8.} {\it The Jordan algebra $M_2({\Bbb C})$ is an RJ-algebra.}

{\it Proof.} By corollary 1.7 in order to prove this proposition it is sufficient to show that
the algebra $M_2({\Bbb C})$ does not have
nilpotent elements, which have square roots.

First we prove that if
a matrix $M$ from $M_2({\Bbb C})$ is nilpotent then $M^2=0$.
Let $M=\left(
  \begin{array}{cc}
    a & b \\
    c & d \\
  \end{array}
\right)$.
Suppose that $M^k=0$ for an integer number $k>2$.
Then
$$
M^2=(a+d)M+(bc-ad)E,\,\,\,\,\,\,\, (1)
$$
where $E$ is the unite matrix.
Since $M$ is nilpotent we have $ad-bc=0$. By (1)
we have the following equality
$$
M^2=(a+d)M.\,\,\,\,\,\,\, (2)
$$
This last equality (2) yields
$$
M^k=(a+d)^{k-1}M.
$$
The equality $M^k=0$ implies $a+d=0$. So $M^2=0$.

Thus, $M^2=0$ for every nilpotent matrix $M$ from $M_2({\Bbb C})$.
Therefore it is sufficient to prove that if $M^2=0$ for a non-zero matrix $M$, then $M$ does not have a square
root. But if $M$ has a square root, say $M=N^2$, then $N^4 = 0$. From the above this implies
that $N^2=0$, i.e. $M=0$  - a contradiction.

Thus the algebra $M_2({\Bbb C})$ does not have
nilpotent elements, which have square roots and $M_2({\Bbb C})$ is an RJ-algebra.
$\triangleright$

\bigskip

\section{BJ-algebras}

\medskip

Let $A$ be a Jordan algebra. Fix the following condition:

(A) for every subset $S\subseteq A$ there exists an idempotent
$e\in A$ such that $^\perp S\cap A^2=U_e(A)\cap A^2$.

This condition is a Jordan analogue of Baer condition for Baer
$*$-algebras.

{\it Definition.} A Jordan algebra satisfying the condition (A) is
called {\it a BJ-algebra}.

If for every subset $S\subseteq A^2$ of a Jordan algebra $A$ there exists an idempotent $e\in
A$ such that $S^\perp =U_e(A)$, then $A$ is called a Baer Jordan algebra \cite{AA}.

{\bf Proposition 2.1.}
{\it Every Baer Jordan algebra is a BJ-algebra.}

{\it Proof.} This proposition follows from theorem 2.6 in \cite{AA}. $\triangleright$

{\it Example.} The Jordan algebras $H_n(F)$ of $n\times n$ dimensional self-adjoint matrices over $F$,
where $F={\Bbb R}$ (reals), ${\Bbb C}$ (complexes), ${\Bbb H}$ (quaternions) or
${\Bbb Q}$ (octonions), are  BJ-algebras by corollary 3.6 in \cite{AA} and proposition 2.1.
By proposition 1.2 and proposition 3.4 in \cite{AA} $H_n(F)$ is an RJ-algebra.
So the set $\mathcal{R}$ of all infinite sequences
with components from $H_n(F)$ is also a BJ-algebra with respect to
componentwise algebraic operations.

Note that the RJ-algebra $\mathcal{R}$ of all infinite
sequences with all but finite non-zero components from $H_n(F)$
is not a BJ-algebra. Because the lattice of
all idempotents of this RJ-algebra is not complete (see lemma 2.6).

{\it Definition.} A Baer $*$-algebra is a $*$-algebra $A$ such
that for every nonempty subset $S$ of $A$, its right annihilator $R(S)=gA$ for a
suitable projection $g$.

The following proposition holds:

{\bf Proposition 2.2.}
{\it Let $A$ be a Baer $*$-algebra. Then
$A_{sa}$ is a BJ-algebra.}

{\bf Proof.} The Jordan algebra $A_{sa}$ is a BJ-algebra
by proposition 2.1 in \cite{AA}. $\triangleright$

{\bf Lemma 2.3.}
{\it Let $A$ be a BJ-algebra. Then

a) $A$ has an element $e$ such that $ae=ea=a$ for every $a\in A^2$,

b) the algebra $A$ has no nilpotent elements, which have square roots. }

{\bf Proof.} This lemma follows from lemma 1.3. $\triangleright$

{\bf Lemma 2.4.}
{\it The set of all idempotents of a BJ-algebra $A$ forms
a complete lattice.}

{\it Proof.}  Let $\{e_i\}_i$ be a set of idempotents in $A$. We have
$$
^\perp \{e_i\}_i\cap A^2=U_{1-e}(A)\cap A^2
$$
for some idempotent $e$ in $A$. Since $^\perp \{e_i\}\cap A^2=U_{1-e_i}(A)\cap A^2$ by lemma 1.4
we have
$$
U_{1-e}(A)\subseteq U_{1-e_i}(A)
$$
for every $i$. Hence $1-e\leq 1-e_i$, i.e. $e_i\leq e$ for each $i$.

We prove that $\sup_i e_i=e$. Indeed,
suppose $f\geq e_i$ for each $i$, where $f$ is an idempotent in $A$.
Then $U_{e_i}(1-f)=0$ for each $i$, i.e. $1-f\in ^\perp \{e_i\}_i\cap A^2=U_{1-e}(A)\cap A^2$.
So $1-f\leq 1-e$, i.e. $e\leq f$. Therefore $\sup_i e_i=e$.

Also we have
$$
\inf_i e_i=1-\sup_i (1-e_i).
$$
So the proof is complete. $\triangleright$

{\bf Theorem 2.5.}
{\it The following conditions are equivalent:

(a) $A$ is a BJ-algebra;

(b) $A$ is an RJ-algebra and the set
of all idempotents of $A$ is a complete lattice.}

{\it Proof.}  (a)$\Longrightarrow$(b): By lemma 2.4 all idempotents of a BJ-algebra
$A$ form a complete lattice. The rest part is obvious.

(b)$\Longrightarrow$(a):  Let $S=\{x_\alpha :\alpha\in \Omega\}$
be any subset of $A$; we have to show that there exists an
idempotent $e\in A$ such that ${^\perp S}\cap A^2=U_e(A)\cap A^2$. Now,
$$
{^\perp S}\cap A^2=\cap_{\alpha\in\Omega} {^\perp\{x_\alpha\}}\cap A^2.
$$
Since $A$ is an RJ-algebra, there exists an idempotent
$e_\alpha\in A$ such that ${^\perp\{x_\alpha\}}\cap A^2=U_{e_\alpha}(A)\cap A^2$ for
any $\alpha$. Let $e=\inf_\alpha e_\alpha$ in $A$, which exists by
the hypothesis of the theorem. Then
$$
e=1-\sup_\alpha(1-e_\alpha)
$$
in $A$. We have
$$
U_e(A)\cap A^2\subseteq {^\perp\{x_\alpha\}}\cap A^2
$$
for any $\alpha$. Hence $U_e(A)\cap A^2\subseteq {^\perp S}\cap A^2$. We prove that
$\cap_{\alpha\in\Omega}{^\perp\{x_\alpha\}}\cap A^2\subseteq U_e(A)\cap A^2$. Let
$a$ be an element in
$\cap_{\alpha\in\Omega}{^\perp\{x_\alpha\}}\cap A^2$. Then $a\in
U_{e_\alpha}(A)\cap A^2$ and $ae_\alpha=a$, i.e. $U_{1-e_\alpha}a=0$ for
any $\alpha$. Hence
$$
a\in {^\perp \{1-e_\alpha\}_\alpha}\cap A^2=U_e(A)\cap A^2, U_{\sup_\alpha (1-e_\alpha)}a=U_{1-e}a=0,
$$
by lemma 1.4 and by the proof of lemma 2.4. Therefore
$a\in {^\perp\{1-e\}}=U_e(A)\cap A^2$ by lemma 1.4, i.e.
$$
\cap_{\alpha\in\Omega}{^\perp\{x_\alpha\}}\cap A^2=^\perp S\cap A^2\subseteq U_e(A)\cap A^2.
$$

Thus ${^\perp S}\cap A^2=U_e(A)\cap A^2$ and $A$ is a BJ-algebra. The proof is complete.
$\triangleright$

{\it Examples.}
1. The RJ-algebra $\mathcal{A}={\Bbb R}{\bf 1}+{\Bbb R}e_{1,2}$
with the multiplication $x\cdot y=\frac{1}{2}(xy+yx)$ is a BJ-algebra by theorem 2.5 since
its lattice of all idempotents $\{1, 0\}$ is complete.
Note that $\mathcal{A}$ is an example of a BJ-algebra which has a
trivial element and a nonzero nilpotent radical.

2. The RJ-algebra
$$
\mathcal{A}={\Bbb R}{\bf 1}+\sum_{i=1}^k {\Bbb R}e_{2i-1,2i}
$$
is also a BJ-algebra. The Jordan algebra $\mathcal{A}$ is also an example of
a BJ-algebra which has a trivial element and a nonzero nilpotent radical.

3. The RJ-algebra $\mathcal{A}=\sum_{i=1}^k {\Bbb R}e_{2i-1,2i}$ without a unit element
is a BJ-algebra by theorem 2.5 since its lattice of all idempotents $\{0\}$ is complete.

{\bf Theorem 2.6.}
{\it Let $A$ be a Jordan algebra. Then, if the quotient
$A/\mathcal{D}eg(A)$ is a BJ-algebra, then $A$ is a BJ-algebra.}

{\it Proof.}
Clearly $A/\mathcal{D}eg(A)$ is an RJ-algebra. Hence by theorem 1.5 $A$ is an RJ-algebra.
By theorem 2.4 the set of all idempotents in $A/\mathcal{D}eg(A)$ forms a complete lattice.
Obviously this lattice is also the set of all idempotents in $A$. Hence by theorem 2.5
$A$ is a BJ-algebra. The proof is complete. $\triangleright$

{\bf Corollary 2.7.} {\it Every finite dimensional Jordan algebra without nilpotent elements
which have square roots is a BJ-algebra.}

{\it Proof.} Let $A$ be a finite dimensional unital Jordan algebra.
By radical equality theorem 1.7.2 in \cite{KMcC} $\mathcal{D}eg(A)=\mathcal{R}ad(J)$.
Hence by enlightenment structure theorem in page 79 of \cite{KMcC}, by theorems 2.4, 3.7 in \cite{AA} and
by proposition 1.2 the quotient $A/\mathcal{D}eg(A)$ is an RJ-algebra. Therefore
by theorem 2.6 $A$ is an RJ-algebra. The proof is complete. $\triangleright$.


\begin{thebibliography}{CS79}


\bibitem{AA} {\em Ayupov Sh.A., Arzikulov F.N.} Jordan counterparts of Rickart and Baer *-algebras.
Uzbek Mathematical Journal. 2016. No 1, 13--33.

\bibitem{JNW} {\em Jordan P., von Neumann J. and Wigner E.} On an algebraic
generalization of the quantum mechanical formalism. Annals
Math. 1934. {\bf 35}. 29--64.

\bibitem{Arz} {\em Arzikulov F.N.} On abstract JW-algebras. Sib. Math.
J. 1998. {\bf 39}. 20--27.

\bibitem{KMP} {\em Kaup W., McCrimmon K., Petersson H.P.} Jordan Algebras: Proceedings of the Conference Held in
Oberwolfach, (Germany, August 9-15, 1992) Walter de Gruyter,
Berlin, New York 1994.

\bibitem{HOS} {\em Hanñhe-Olsen H., St$\hat{o}$rmer E.} Jordan Operator Algebras. Boston
etc.:Pitman Publ. Inc. 1984.

\bibitem{KMcC} {\em McCrimmon K.} A test of Jordan algebras.
Springer, New York and etc, 2004, P. 562.


\end{thebibliography}
\end{document}